\documentclass{tac}

\usepackage{amsmath}
\usepackage{amssymb}
\usepackage{enumitem}

\usepackage[all,cmtip]{xy}

\input diagxy

\usepackage[colorlinks=true]{hyperref}
\hypersetup{allcolors=[rgb]{0.1,0.1,0.4}}

\author{Lili Shen and Hang Yang}

\thanks{The authors acknowledge the support of National Natural Science Foundation of China (No. 12071319).}

\address{School of Mathematics, Sichuan University\\
 Chengdu 610064, China\\
}

\title{Injective symmetric quantaloid-enriched categories}

\copyrightyear{2022}


\eaddress{shenlili@scu.edu.cn\CR yanghangscu@qq.com}

\newtheorem{thm}{Theorem}

\newtheorem{lem}{Lemma}
\newtheorem{prop}{Proposition}
\newtheorem{cor}{Corollary}

\newtheoremrm{rem}{Remark}
\newtheoremrm{defn}{Definition}
\newtheoremrm{exmp}{Example}
\newtheoremrm{ques}{Question}

\mathrmdef{Hom}

\DeclareMathOperator{\id}{id}
\DeclareMathOperator{\ob}{ob}

\def\oto{{\bfig\morphism<180,0>[\mkern-4mu`\mkern-4mu;]\place(86,0)[\circ]\efig}}
\def\rto{{\bfig\morphism<180,0>[\mkern-4mu`\mkern-4mu;]\place(82,0)[\mapstochar]\efig}}

\newcommand{\ra}{\rightarrow}

\newcommand{\lda}{\swarrow}
\newcommand{\rda}{\searrow}
\newcommand{\bs}{\backslash}

\newcommand{\bv}{\bigvee}
\newcommand{\bw}{\bigwedge}

\newcommand{\dv}{\dashv}

\newcommand{\nat}{\natural}

\renewcommand{\phi}{\varphi}
\newcommand{\al}{\alpha}
\newcommand{\be}{\beta}

\newcommand{\lam}{\lambda}
\newcommand{\si}{\sigma}

\newcommand{\Om}{\Omega}

\newcommand{\CQ}{\mathcal{Q}}

\newcommand{\CW}{\mathcal{W}}

\newcommand{\sD}{{\sf D}}

\newcommand{\sL}{{\sf L}}

\newcommand{\sP}{{\sf P}}
\newcommand{\sQ}{{\sf Q}}

\newcommand{\sT}{{\sf T}}

\newcommand{\sy}{{\sf y}}
\newcommand{\sfs}{{\sf s}}

\newcommand{\Cat}{{\bf Cat}}

\newcommand{\Dist}{{\bf Dist}}

\newcommand{\Met}{{\bf Met}}

\newcommand{\Rel}{{\bf Rel}}
\newcommand{\Set}{{\bf Set}}
\newcommand{\Sup}{{\bf Sup}}

\newcommand{\ParMet}{{\bf ParMet}}

\newcommand{\QCat}{\CQ\text{-}\Cat}
\newcommand{\QSymCat}{\CQ\text{-}{\bf SymCat}}

\newcommand{\QDist}{\CQ\text{-}\Dist}

\newcommand{\QRel}{\CQ\text{-}\Rel}

\newcommand{\sPd}{\sP^{\dag}}

\newcommand{\op}{{\rm op}}

\newcommand{\PX}{\sP X}

\newcommand{\LX}{\sL X}
\newcommand{\LY}{\sL Y}

\newcommand{\TX}{\sT X}
\newcommand{\TY}{\sT Y}

\newcommand{\PdX}{\sPd X}

\newcommand{\DQ}{\sD\sQ}

\newcommand{\ldd}{\mathrel{/}}
\newcommand{\rdd}{\mathrel{\bs}}

\renewcommand{\leq}{\leqslant}
\renewcommand{\geq}{\geqslant}

\newcommand{\tmu}{\widetilde{\mu}}
\newcommand{\tlam}{\widetilde{\lam}}

\numberwithin{equation}{section}

\begin{document}

\maketitle
\begin{abstract}
We characterize injective objects, injective hulls and essential embeddings in the category of symmetric categories enriched in a small, integral and involutive quantaloid. In particular, injective partial metric spaces are precisely formulated.
\end{abstract}


\section{Introduction}

A \emph{quantaloid} \cite{Rosenthal1996} $\CQ$ is a category enriched in the monoidal-closed category $\Sup$ \cite{Joyal1984} of complete lattices and $\sup$-preserving maps. Considering $\CQ$ as a base for enrichment, a theory of $\CQ$-categories, $\CQ$-functors and $\CQ$-distributors can be developed \cite{Rosenthal1996,Stubbe2005,Stubbe2006}. Moreover, if $\CQ$ is involutive, it makes sense to consider \emph{symmetric} $\CQ$-categories \cite{Betti1982a,Heymans2011}. 

This paper is concerned with \emph{injectivity} \cite{Maranda1964,Adamek1990} in the category of symmetric $\CQ$-categories. More specifically, since the Lawvere quantale \cite{Lawvere1973}
$$[0,\infty]=([0,\infty],+,0)$$
is trivially an involutive quantaloid, and (classical) metric spaces are symmetric $[0,\infty]$-categories, we wish to find the categorical interpretation of \emph{injective metric spaces} in the framework of quantaloid-enriched categories.

It is well known that injective metric spaces, i.e., injective objects in the category $\Met$ of (classical) metric spaces and non-expansive maps, are precisely \emph{hyperconvex} metric spaces:

\begin{thm} \label{Met-injective} {\rm\cite{Aronszajn1956,Isbell1964,Adamek1990,Espinola2001}}
A metric space $(X,\al)$ is injective if, and only if, it is hyperconvex in the sense that for every family $\{(x_j,r_j)\}_{j\in J}$ of pairs with $x_j\in X$ and $r_j\in[0,\infty]$ satisfying
$$\al(x_j,x_k)\leq r_j+r_k$$
for all $j,k\in J$, there exists $z\in X$ such that
$$\al(x_j,z)\leq r_j$$
for all $j\in J$.
\end{thm}

Furthermore, the \emph{injective hull} (also \emph{injective envelope}) of a metric space was firstly constructed by Isbell \cite{Isbell1964}, and later characterized by Dress \cite{Dress1984} as the \emph{tight span}:

\begin{thm} \label{Met-injective-hull} {\rm\cite{Isbell1964,Dress1984,Herrlich1992a,Willerton2013}}
The injective hull of a metric space $(X,\al)$ is given by its tight span $(\sT(X,\al),\si)$, where $\sT(X,\al)$ consists of maps $\mu:X\to[0,\infty]$ satisfying
$$\mu(x)=\sup\limits_{y\in X}(\al(x,y)-\mu(y))$$
for all $x\in X$, and
$$\si(\mu,\lam)=\sup\limits_{x\in X}(\mu(x)-\lam(x))$$
for all $\mu,\lam\in\sT(X,\al)$.
\end{thm}

In this paper we look deeply into the categorical meaning of the concepts of hyperconvexity and tight span, and provide a far reaching extension of Theorems \ref{Met-injective} and \ref{Met-injective-hull} in the context of quantaloid-enriched categories. Explicitly, for a small, integral and involutive quantaloid $\CQ$, where the integrality means that the top element of each $\CQ(q,q)$ is $1_q$, we explore injectivity with respect to \emph{fully faithful $\CQ$-functors} (considered as ``embeddings'') in the category
$$\QSymCat$$
of symmetric $\CQ$-categories and $\CQ$-functors. The main results of this paper are as follows:
\begin{itemize}
\item We define \emph{hyperconvexity} for a symmetric $\CQ$-category (Definition \ref{hyperconvex-def}), and show that a symmetric $\CQ$-category is hyperconvex if, and only if, it is an injective object in $\QSymCat$ (Theorem \ref{hyperconvex=injective}).
\item We define the \emph{tight span} of a symmetric $\CQ$-category (Equations \eqref{LX-TX-def}), and show that it is precisely the injective hull of a symmetric $\CQ$-category in $\QSymCat$ (Theorem \ref{injective-hull}).
\item We show that essential embeddings in $\QSymCat$ are precisely (co)dense fully faithful $\CQ$-functors (Theorem \ref{f-essential}).
\end{itemize}
Finally, in Section \ref{Example}, we demonstrate that our results can be applied to a large family of small, integral and involutive quantaloids constructed upon integral and involutive quantales $\sQ$; that is, the quantaloid $\DQ$ of \emph{diagonals} of $\sQ$ \cite{Hoehle2011a,Pu2012,Stubbe2014}. In particular, we postulate the concepts of \emph{hyperconvexity} and \emph{tight span} for partial metric spaces of Matthews \cite{Matthews1994} (Definition \ref{parmet-hyperconvex-def}, Corollaries \ref{parmet-hyperconvex=injective} and \ref{tight-span}).

It is noteworthy to point out that, although some techniques used in the classical proofs of Theorems \ref{Met-injective} and \ref{Met-injective-hull} are inevitable in our general setting (mostly regarding the use of Zorn's lemma to prove the existence of certain elements), our methods are quite different from that special case. While the classical proofs of Theorems \ref{Met-injective} and \ref{Met-injective-hull} were mostly formulated in a geometric way, in this paper we avoid pointwise computations and always proceed by ``distributor calculus'', which has the potential for a wider range of applicability.

\section{Symmetric quantaloid-enriched categories} \label{Sym-Q-Cat}

A \emph{quantaloid} \cite{Rosenthal1996} is a category enriched in the monoidal-closed category $\Sup$ \cite{Joyal1984} of complete lattices and $\sup$-preserving maps. Explicitly, a quantaloid $\CQ$ is a category whose hom-sets are complete lattices, such that the composition of morphisms preserves arbitrary suprema on both sides. The corresponding right adjoints induced by the compositions 
$$-\circ u\dv -\lda u:\ \CQ(p,r)\to\CQ(q,r)\quad\text{and}\quad v\circ -\dv v\rda -:\ \CQ(p,r)\to\CQ(p,q)$$ 
satisfy 
$$v\circ u\leq w\iff v\leq w\lda u\iff u\leq v\rda w$$ 
for all morphisms $u:p\to q$, $v:q\to r$, $w:p\to r$ in $\CQ$, where $\lda$ and $\rda$ are called \emph{left} and \emph{right implications} in $\CQ$, respectively.

A \emph{homomorphism} of quantaloids is a functor of the underlying categories that preserves suprema of morphisms. A quantaloid $\CQ$ is \emph{involutive} if it is equipped with an \emph{involution}; that is, a homomorphism
$$(-)^{\circ}:\CQ^{\op}\to\CQ$$
of quantaloids with
$$q^{\circ}=q\quad\text{and}\quad u^{\circ\circ}=u$$
for all $q\in\ob\CQ$ and morphisms $u:p\to q$ in $\CQ$, which necessarily satisfies
\begin{equation} \label{lda-rda-circ}
(w\lda u)^{\circ}=u^{\circ}\rda w^{\circ}\quad\text{and}\quad(v\rda w)^{\circ}=w^{\circ}\lda v^{\circ}
\end{equation}
for all morphisms $u,u_j:p\to q$ $(j\in J)$, $v:q\to r$, $w:p\to r$ in $\CQ$.

Throughout this paper, we let $\CQ$ denote a \emph{small}, \emph{integral} and involutive quantaloid, where the integrality of $\CQ$ means that the identity morphism $1_q$ is the top element of the complete lattice $\CQ(q,q)$ for all $q\in\ob\CQ$.

From $\CQ$ we form a (large) involutive quantaloid $\QRel$ of \emph{$\CQ$-relations} with the following data:
\begin{itemize}
\item objects of $\QRel$ are those of $\Set/\ob\CQ$, i.e., sets $X$ equipped with a \emph{type} map $|\text{-}|:X\to\ob\CQ$;
\item a morphism $\phi:X\rto Y$ in $\QRel$ is a family of morphisms $\phi(x,y):|x|\to|y|$ $(x,y\in X)$ in $\CQ$, and its composite with $\psi:Y\rto Z$ is given by
$$\psi\circ\phi:X\rto Z,\quad(\psi\circ\phi)(x,z)=\bv_{y\in Y}\psi(y,z)\circ\phi(x,y),$$
with 
$$\id_X:X\rto X,\quad \id_X(x,y)=\begin{cases}
1_{|x|} & \text{if}\ x=y,\\
\bot & \text{else}
\end{cases}$$
serving as the identity morphism on $X$;
\item the local order in $\QRel$ is inherited from $\CQ$, i.e.,
$$\phi\leq\phi':X\rto Y\iff\forall x\in X,\ \forall y\in Y:\ \phi(x,y)\leq\phi'(x,y);$$
\item implications in $\QRel$ are computed pointwise as
\begin{align*}
&\xi\lda\phi:Y\rto Z,\quad(\xi\lda\phi)(y,z)=\bw_{x\in X}\xi(x,z)\lda\phi(x,y),\\
&\psi\rda\xi:X\rto Y,\quad(\psi\rda\xi)(x,y)=\bw_{z\in Z}\psi(y,z)\rda\xi(x,z)
\end{align*}
for all $\phi:X\rto Y$, $\psi:Y\rto Z$, $\xi:X\rto Z$ in $\QRel$.
\item the involution on $\QRel$ is given by 
$$\phi^{\circ}:Y\rto X,\quad\phi^{\circ}(y,x):=\phi(x,y)^{\circ}$$
for all $\phi:X\rto Y$ in $\QRel$.
\end{itemize}

A \emph{$\CQ$-category} \cite{Rosenthal1996,Stubbe2005} is an internal monad in $\QRel$; that is, an object in $\Set/\ob\CQ$ equipped with a $\CQ$-relation $\al:X\rto X$ with $\id_X\leq\al$ and $\al\circ\al\leq\al$. For every $\CQ$-category $(X,\al)$, the underlying (pre)order on $X$ is given by
$$x\leq y\iff|x|=|y|\ \ \text{and}\ \ \al(x,y)=1_{|x|},$$
and we write $x\cong y$ if $x\leq y$ and $y\leq x$.  It is clear that
\begin{equation} \label{x-cong-y-al-x-y}
x\cong y\iff\al(x,-)=\al(y,-)\iff\al(-,x)=\al(-,y)
\end{equation}
for all $x,y\in X$. We say that $(X,\al)$ is \emph{separated} (also \emph{skeletal}) if $x=y$ whenever $x\cong y$ in $X$.

A map $f:(X,\al)\to(Y,\be)$ between $\CQ$-categories is a \emph{$\CQ$-functor} if
$$|x|=|fx|\quad\text{and}\quad\al(x,x')\leq\be(fx,fx')$$
for all $x,x'\in X$, and $f$ is \emph{fully faithful} if $\al(x,x')=\be(fx,fx')$ for all $x,x'\in X$. With the pointwise order of $\CQ$-functors inherited from $Y$, i.e.,
$$f\leq g:(X,\al)\to(Y,\be)\iff\forall x\in X:\ fx\leq gx\iff\forall x\in X:\ \be(fx,gx)=1_{|x|},$$
$\CQ$-categories and $\CQ$-functors are organized into a 2-category $\QCat$.

A $\CQ$-category $(X,\al)$ is \emph{symmetric} \cite{Heymans2011} if $\al^{\circ}=\al$, i.e.,
\begin{equation} \label{sym-Q-cat}
\al(x,y)=\al(y,x)^{\circ}
\end{equation}
for all $x,y\in X$. The full subcategory of $\QCat$ consisting of symmetric $\CQ$-categories is denoted by 
$$\QSymCat.$$
It is well known that $\QSymCat$ is a coreflective subcategory of $\QCat$ \cite{Heymans2011}, with the coreflector sending each $\CQ$-category $(X,\al)$ to its \emph{symmetrization} $(X,\al_{\sfs})$, where $\al_{\sfs}=\al\wedge\al^{\circ}$, i.e.,
\begin{equation} \label{symmetrization-def}
\al_{\sfs}(x,y)=\al(x,y)\wedge\al(y,x)^{\circ}
\end{equation}
for all $x,y\in X$. Indeed, the coreflectivity of $\QSymCat$ in $\QCat$ is easily observed from the following lemma:

\begin{lem} \label{QSymCat-coref-QCat}
Let $(X,\al)$ be a symmetric $\CQ$-category, and let $(Y,\be)$ be a $\CQ$-category. Then $f:(X,\al)\to(Y,\be)$ is a $\CQ$-functor if, and only if, $f:(X,\al)\to(Y,\be_{\sfs})$ is a $\CQ$-functor.
\end{lem}

A $\CQ$-relation $\phi:X\rto Y$ becomes a \emph{$\CQ$-distributor} $\phi:(X,\al)\oto(Y,\be)$ if
$$\be\circ\phi\circ\al\leq\phi.$$
$\CQ$-categories and $\CQ$-distributors constitute a (\emph{not} necessarily involutive!) quantaloid $\QDist$ that contains $\QRel$ as a full subquantaloid, in which the compositions and implications are calculated in the same way as in $\QRel$, and the identity $\CQ$-distributor on $(X,\al)$ is given by $\al:(X,\al)\oto(X,\al)$.

Each $\CQ$-functor $f:(X,\al)\to(Y,\be)$ induces an adjunction $f_{\nat}\dv f^{\nat}$ in $\QDist$ (i.e., $\al\leq f^{\nat}\circ f_{\nat}$ and $f_{\nat}\circ f^{\nat}\leq\be$) given by
\begin{align*}
&f_{\nat}:(X,\al)\oto(Y,\be),\quad f_{\nat}(x,y)=\be(fx,y),\\
&f^{\nat}:(Y,\be)\oto(X,\al),\quad f^{\nat}(y,x)=\be(y,fx),
\end{align*}
called the \emph{graph} and \emph{cograph} of $f$, respectively. Obviously, the identity $\CQ$-distributor $\al:(X,\al)\oto(X,\al)$ is the cograph of the identity $\CQ$-functor $1_X:X\to X$. Hence, in what follows
$$1_X^{\nat}=\al$$
will be our standard notation for the hom of a $\CQ$-category $X=(X,\al)$ if no confusion arises. 

It is straightforward to verify the following lemmas:

\begin{lem} \label{fully-faithful-graph} {\rm\cite{Shen2013a}}
A $\CQ$-functor $f:X\to Y$ is fully faithful if, and only if, $f^{\nat}\circ f_{\nat}=1_X^{\nat}$.
\end{lem}

\begin{lem} \label{graph-cograph-calc} 
The following identities hold for all $\CQ$-functors $f$, $g$ and $\CQ$-distributors $\phi$, $\psi$ whenever the operations make sense:
\begin{enumerate}[label={\rm(\arabic*)}]
\item \label{graph-cograph-calc:composition} $(gf)_{\nat}=g_{\nat}\circ f_{\nat}$,\quad $(gf)^{\nat}=f^{\nat}\circ g^{\nat}$.
\item \label{graph-cograph-calc:phi-f-g} $\phi(f-,g-)=g^{\nat}\circ\phi\circ f_{\nat}$.
\item {\rm\cite{Heymans2010}} \label{graph-cograph-calc:f-nat-bracket}  $(\phi\rda\psi)\circ f_{\nat}=\phi\rda(\psi\circ f_{\nat})$,\quad $f^{\nat}\circ(\psi\lda\phi)=(f^{\nat}\circ\psi)\lda\phi$.
\end{enumerate}
\end{lem}

For each $q\in\ob\CQ$, let $\{q\}$ denote the (necessarily symmetric) one-object $\CQ$-category whose only object has type $q$ and hom $1_q$. A \emph{presheaf} $\mu$ (of type $q$) on a $\CQ$-category $X$ is a $\CQ$-distributor $\mu:X\oto\{q\}$, and presheaves on $X$ constitute a separated $\CQ$-category $\PX$ with
$$1_{\PX}^{\nat}(\mu,\mu')=\mu'\lda\mu$$
for all $\mu,\mu'\in\PX$. Dually, the separated $\CQ$-category $\PdX$ of \emph{copresheaves} on $X$ consists of $\CQ$-distributors $\lam:\{q\}\oto X$ with $|\lam|=q$ and
$$1_{\PdX}^{\nat}(\lam,\lam')=\lam'\rda\lam$$
for all $\lam,\lam'\in\PdX$.

For each $\CQ$-distributor $\phi:X\oto Y$ between symmetric $\CQ$-categories, it is easy to see that
$$\phi^{\circ}:Y\oto X,\quad\phi^{\circ}(y,x):=\phi(x,y)^{\circ}$$
is also a $\CQ$-distributor. Therefore:
\begin{itemize}
\item For a symmetric $\CQ$-category $X$, $\mu:X\oto\{q\}$ is a presheaf on $X$ if, and only if, $\mu^{\circ}:\{q\}\oto X$ is a copresheaf on $X$. 
\item For a $\CQ$-functor $f:X\to Y$ between symmetric $\CQ$-categories,
\begin{equation} \label{f-graph-circ}
f_{\nat}^{\circ}=f^{\nat}\quad\text{and}\quad(f^{\nat})^{\circ}=f_{\nat}.
\end{equation}
\end{itemize}

\section{Injective symmetric quantaloid-enriched categories} \label{Inj-Sym-Q-Cat}




While discussing the injectivity of $\CQ$-categories, it is standard to consider injective objects with respect to fully faithful $\CQ$-functors, which may be treated as ``embeddings'' of $\CQ$-categories; we refer to \cite{Stubbe2006a,Hofmann2011b,Stubbe2017,Shen2016,Fujii2021} for the counterparts in $\QCat$ of our main results:
\begin{itemize}
\item A $\CQ$-category is an injective object in $\QCat$ if, and only if, it is \emph{(co)complete}.
\item The injective hull of a $\CQ$-category in $\QCat$ is its \emph{MacNeille completion}.
\item Essential embeddings in $\QCat$ are precisely \emph{dense} and \emph{codense} fully faithful $\CQ$-functors.
\end{itemize}

The aim of this section is to characterize injective objects (with respect to fully faithful $\CQ$-functors) in the category $\QSymCat$. To clarify the terminology, we say that a symmetric $\CQ$-category $Z$ is \emph{injective} in $\QSymCat$ if, for every $\CQ$-functor $f:X\to Z$ and every fully faithful $\CQ$-functor $g:X\to Y$ between symmetric $\CQ$-categories, there exists a (not necessarily unique) $\CQ$-functor $h:Y\to Z$ such that $hg\cong f$.
$$\bfig
\qtriangle/->`->`-->/<600,400>[X`Y`Z;g`f`h]
\efig$$

\begin{defn} \label{hyperconvex-def}
A symmetric $\CQ$-category $X$ is \emph{hyperconvex} if, for every presheaf $\mu$ on $X$ satisfying
\begin{equation} \label{mu-mu-leq-1X}
\mu^{\circ}\circ\mu\leq 1_X^{\nat},
\end{equation}
there exists $z\in X$ such that
$$\mu\leq 1_X^{\nat}(-,z).$$
\end{defn}

\begin{rem}
The reason that we use the term ``hyperconvex'' is obvious. In the case that $\CQ=[0,\infty]$ is the Lawvere quantale (see Example \ref{quantale-exmp}\ref{quantale-exmp:Lawvere}), a (classical) metric space is \emph{hyperconvex} if, and only if, it is hyperconvex as a symmetric $[0,\infty]$-category. 
\end{rem}

\begin{rem} \label{hyperconvex-def-QRel}
For the sake of defining hyperconvex partial metric spaces in Section \ref{Example}, we point out that a symmetric $\CQ$-category $X$ is hyperconvex if, and only if, for every $\CQ$-relation $\mu:X\rto\{q\}$ (not necessarily a presheaf on $X$) satisfying
$$\mu^{\circ}\circ\mu\leq 1_X^{\nat},$$
there exists $z\in X$ such that
$$\mu\leq 1_X^{\nat}(-,z).$$
Indeed, every $\CQ$-relation $\mu:X\rto\{q\}$ gives rise to a presheaf $\mu\circ 1_X^{\nat}$ on $X$, and $\mu$ satisfies any of the above inequalities if, and only if, so does $\mu\circ 1_X^{\nat}$.
\end{rem}

Let $Y$ be a $\CQ$-category. Then each subset $X\subseteq Y$ is equipped with the $\CQ$-categorical structure inherited from $Y$, i.e.,
$$1_X^{\nat}(x,x')=1_Y^{\nat}(x,x')$$
for all $x,x'\in X$. In this case, we say that $X$ is a \emph{$\CQ$-subcategory} of $Y$, and $Y$ is a \emph{$\CQ$-supercategory} of $X$.

A $\CQ$-category $X$ is a \emph{retract} of a $\CQ$-category $Y$ if there are $\CQ$-functors $i:X\to Y$ and $h:Y\to X$ such that $hi\cong 1_X$. In this case, $h$ is called a \emph{retraction} of $Y$ onto $X$, and $i$ is called a \emph{section} of $h$. When $X$ is symmetric, we say that
\begin{itemize}
\item $X$ is a \emph{pre-absolute retract} in $\QSymCat$, if the inclusion $\CQ$-functor $i:X\ \to/^(->/Y$ is a section whenever $Y$ is symmetric $\CQ$-supercategory of $X$ with $Y\setminus X$ containing exactly one element;
\item $X$ is an \emph{absolute retract} in $\QSymCat$, if the inclusion $\CQ$-functor $i:X\ \to/^(->/Y$ is a section whenever $Y$ is symmetric $\CQ$-supercategory of $X$.
\end{itemize}

\begin{prop} \label{retract-hyperconvex}
If $X$ is a retract of a hyperconvex symmetric $\CQ$-category $Y$ in $\QSymCat$, then $X$ is also hyperconvex.
\end{prop}

\begin{proof}
Suppose that there are $\CQ$-functors $i:X\to Y$ and $h:Y\to X$ satisfying $hi\cong 1_X$. Let $\mu$ be a presheaf on $X$ satisfying \eqref{mu-mu-leq-1X}. Then $\mu\circ i^{\nat}$ is a presheaf on $Y$, and it follows from Equations \eqref{f-graph-circ} and $i_{\nat}\dv i^{\nat}$ that
$$(\mu\circ i^{\nat})^{\circ}\circ(\mu\circ i^{\nat})=i_{\nat}\circ \mu^{\circ}\circ\mu\circ i^{\nat}\leq i_{\nat}\circ i^{\nat}\leq 1_Y^{\nat}.$$
Hence, by the hyperconvexity of $Y$, there exists $z\in Y$ such that
$$\mu\circ i^{\nat}\leq 1_Y^{\nat}(-,z),$$
and consequently, the $\CQ$-functoriality of $h$ implies that
$$\mu\leq 1_Y^{\nat}(-,z)\lda i^{\nat}=1_Y^{\nat}(-,z)\lda 1_Y^{\nat}(-,i-)=1_Y^{\nat}(i-,z)\leq 1_X^{\nat}(hi-,hz)=1_X^{\nat}(-,hz),$$
which shows that $X$ is also hyperconvex.
\end{proof}

Without the axiom of choice we may characterize hyperconvex symmetric $\CQ$-categories as follows:

\begin{prop} \label{hyperconvex=one-point-retract}
A symmetric $\CQ$-category $X$ is hyperconvex if, and only if, $X$ is a a pre-absolute retract in $\QSymCat$.
\end{prop}

\begin{proof}
``$\implies$'': Suppose that $X$ is hyperconvex. Let $Y$ be a symmetric $\CQ$-supercategory of $X$ with $Y\setminus X=\{y_0\}$. Let $i:X\ \to/^(->/Y$ be the inclusion $\CQ$-functor. Then the fully faithfulness of $i$ ensures that $i^{\nat}\circ i_{\nat}=1_X^{\nat}$ (Lemma \ref{fully-faithful-graph}), and it follows from Equations \eqref{f-graph-circ} that the presheaf $i_{\nat}(-,y_0)$ on $X$ satisfies
$$(i_{\nat}(-,y_0))^{\circ}\circ i_{\nat}(-,y_0)=i^{\nat}(y_0,-)\circ i_{\nat}(-,y_0)\leq i^{\nat}\circ i_{\nat}=1_X^{\nat}.$$
Hence, the hyperconvexity of $X$ guarantees the existence of $z_0\in X$ such that 
$$i_{\nat}(-,y_0)\leq 1_X^{\nat}(-,z_0).$$ 
Define
\begin{equation} \label{h-def}
h:Y\to X,\quad hy=\begin{cases}
y & \text{if}\ y\in X,\\
z_0 & \text{if}\ y=y_0.
\end{cases}
\end{equation}
Then $h:Y\to X$ is a $\CQ$-functor, because 
\begin{equation} \label{h-functor}
1_Y^{\nat}(y_0,y_0)=1_X^{\nat}(z_0,z_0)=1_X^{\nat}(hy_0,hy_0)=1_{|y_0|}
\end{equation}
and
$$1_Y^{\nat}(x,y_0)=i_{\nat}(x,y_0)\leq 1_X^{\nat}(x,z_0)=1_X^{\nat}(hx,hy_0)$$
for all $x\in X$. From the definition of $h$ we immediately deduce that $h|_X=1_X$ ($h|_X$ refers to the restriction of $h$ on $X$); that is, the inclusion $\CQ$-functor $i:X\ \to/^(->/Y$ is a section.

``$\impliedby$'': Let $\mu$ be a presheaf on $X$ satisfying 
$$\mu^{\circ}\circ\mu\leq 1_X^{\nat}.$$
Define a $\CQ$-supercategory 
$$Y:=X\cup\{\mu\}$$ 
of $X$ with
$$1_Y^{\nat}(x,\mu)=\mu(x),\quad 1_Y^{\nat}(\mu,x)=\mu^{\circ}(x),\quad 1_Y^{\nat}(\mu,\mu)=1_{|\mu|}$$
for all $x\in X$. Then $Y$ is symmetric, and the inclusion $\CQ$-functor $i:X\ \to/^(->/Y$ is a section; that is, there exists a $\CQ$-functor $h:Y\to X$ such that $h|_X\cong 1_X$. We claim that 
$$\mu\leq 1_X^{\nat}(-,h\mu).$$ 
Indeed, the $\CQ$-functoriality of $h$ implies that
$$\mu(x)=1_Y^{\nat}(x,\mu)\leq 1_X^{\nat}(hx,h\mu)=1_X^{\nat}(x,h\mu)$$
for all $x\in X$, which completes the proof.
\end{proof}



Assuming the axiom of choice, we arrive at the main result of this section:

\begin{thm} \label{hyperconvex=injective}
Let $X$ be a symmetric $\CQ$-category. Then the following statements are equivalent:
\begin{enumerate}[label={\rm(\roman*)}]
\item \label{hyperconvex=injective:hyperconvex} $X$ is hyperconvex.
\item \label{hyperconvex=injective:one-point} $X$ is a pre-absolute retract in $\QSymCat$.
\item \label{hyperconvex=injective:AR} $X$ is an absolute retract in $\QSymCat$.
\item \label{hyperconvex=injective:injective} $X$ is an injective object in $\QSymCat$.
\end{enumerate}
\end{thm}

\begin{proof}
Since \ref{hyperconvex=injective:injective}$\implies$\ref{hyperconvex=injective:AR}$\implies$\ref{hyperconvex=injective:one-point} is trivial and \ref{hyperconvex=injective:hyperconvex}$\iff$ \ref{hyperconvex=injective:one-point} is already obtained in Proposition \ref{hyperconvex=one-point-retract}, it remains to prove \ref{hyperconvex=injective:hyperconvex}$\implies$\ref{hyperconvex=injective:injective}.

For symmetric $\CQ$-categories $X,Y,Z$, suppose that $Z$ is hyperconvex, $f:X\to Z$ is a $\CQ$-functor and $g:X\to Y$ is a fully faithful $\CQ$-functor. In what follows we use Zorn's lemma to prove the existence of a $\CQ$-functor $h:Y\to Z$ such that $hg\cong f$.

Let $Y_0$ be the $\CQ$-subcategory of $Y$ given by
$$Y_0:=\{gx\mid x\in X\}.$$
Then, for each $y\in Y_0$ we may choose $g^*y\in X$ such that $gg^*y=y$, making $g^*:Y_0\to X$ a fully faithful $\CQ$-functor. It is easy to check that $g^*g\cong 1_X$ through \eqref{x-cong-y-al-x-y}, and consequently $fg^*g\cong f$.

Now, let $\CW$ denote the set of pairs $(W,h_W)$, where $Y_0\subseteq W\subseteq Y$, $h_W:W\to Z$ is a $\CQ$-functor, such that $h_W g\cong f$. We have already seen that $(Y_0,g^*g)\in\CW$, and thus $\CW\neq\varnothing$. Note that $\CW$ is equipped with a partial order given by 
$$(W,h_W)\leq(W',h_{W'})\iff W\subseteq W'\ \text{and}\ h_{W'}|_W=h_W,$$
where $h_{W'}|_W$ refers to the restriction of $h_{W'}$ on $W$. Then, it is easy to see that every chain $\{(W_j,h_{W_j})\mid j\in J\}\subseteq\CW$ has an upper bound in $\CW$, given by 
$$\Big(\bigcup\limits_{j\in J}W_j,h_{\bigcup\limits_{j\in J}W_j}\Big).$$
By Zorn's lemma, $\CW$ has a maximal element $(\overline{W},h_{\overline{W}})$. We claim that $\overline{W}=Y$, which guarantees the existence of the required $h=h_{\overline{W}}:Y\to Z$.

We proceed by contradiction. Assume that there exists $y_0\in Y$ but $y_0\not\in\overline{W}$. Let $W_0=\overline{W}\cup\{y_0\}$. Consider the presheaf $\mu$ on $Z$ given by
$$\bfig
\morphism<600,0>[\mu=(Z`\overline{W};h_{\overline{W}}^{\nat}]
\morphism(600,0)<700,0>[\overline{W}`\{|y_0|\}),;i_{\nat}(-,y_0)]
\place(350,0)[\circ] \place(900,0)[\circ]
\efig$$
where $i:\overline{W}\ \to/^(->/W_0$ is the inclusion $\CQ$-functor. Then
\begin{align*}
\mu^{\circ}\circ\mu&=(h_{\overline{W}}^{\nat})^{\circ}\circ(i_{\nat}(-,y_0))^{\circ}\circ i_{\nat}(-,y_0)\circ h_{\overline{W}}^{\nat}\\
&=(h_{\overline{W}})_{\nat}\circ i^{\nat}(y_0,-)\circ i_{\nat}(-,y_0)\circ h_{\overline{W}}^{\nat}&(\text{Equations \eqref{f-graph-circ}})\\
&\leq(h_{\overline{W}})_{\nat}\circ i^{\nat}\circ i_{\nat}\circ h_{\overline{W}}^{\nat}\\
&=(h_{\overline{W}})_{\nat}\circ h_{\overline{W}}^{\nat}&(\text{Lemma \ref{fully-faithful-graph}})\\
&\leq 1_Z^{\nat}.&((h_{\overline{W}})_{\nat}\dv h_{\overline{W}}^{\nat})
\end{align*}
Hence, the hyperconvexity of $Z$ guarantees the existence of $z_0\in Z$ such that 
$$\mu\leq 1_Z^{\nat}(-,z_0).$$ 
Define
\begin{equation} \label{hW0-def}
h_{W_0}:W_0\to Z,\quad h_{W_0}y=\begin{cases}
h_{\overline{W}}y & \text{if}\ y\in\overline{W},\\
z_0 & \text{if}\ y=y_0.
\end{cases}
\end{equation}
Since $gx\in Y_0\subseteq\overline{W}$ for all $x\in X$ and $h_{\overline{W}}g\cong f$, it is clear that $h_{W_0}g\cong f$. Moreover, $h_{W_0}:W_0\to Z$ is a $\CQ$-functor, because 
\begin{equation} \label{hW0-functor}
1_{W_0}^{\nat}(y_0,y_0)=1_Z^{\nat}(z_0,z_0)=1_Z^{\nat}(h_{W_0}y_0,h_{W_0}y_0)=1_{|y_0|}
\end{equation}
and
\begin{align*}
1_{W_0}^{\nat}(y,y_0)&=i_{\nat}(y,y_0)&(y\in\overline{W}\ \text{and}\ i:\overline{W}\ \to/^(->/W_0)\\
&=i_{\nat}(-,y_0)\circ 1_{\overline{W}}^{\nat}(y,-)\\
&\leq i_{\nat}(-,y_0)\circ 1_Z^{\nat}(h_{\overline{W}}y,h_{\overline{W}}-)&(h_{\overline{W}}\ \text{is a}\ \CQ\text{-functor})\\
&=i_{\nat}(-,y_0)\circ h_{\overline{W}}^{\nat}(h_{\overline{W}}y,-)\\
&=\mu h_{\overline{W}}y&(\text{definition of}\ \mu)\\
&\leq 1_Z^{\nat}(h_{\overline{W}}y,z_0)\\
&=1_Z^{\nat}(h_{W_0}y,h_{W_0}y_0)
\end{align*}
for all $y\in\overline{W}$. Therefore, $(W_0,h_{W_0})\in\CW$, contradicting to the maximality of $(\overline{W},h_{\overline{W}})$, which completes the proof.
\end{proof}

\begin{rem}
In the proofs of Proposition \ref{hyperconvex=one-point-retract} and Theorem \ref{hyperconvex=injective}, the integrality of the quantaloid $\CQ$ is applied to the $\CQ$-functoriality of \eqref{h-def} and \eqref{hW0-def}. In fact, \eqref{h-functor} holds because both $1_Y^{\nat}(y_0,y_0)$ and $1_X^{\nat}(z_0,z_0)$ are \emph{equal} to $1_{|y_0|}$, and \eqref{hW0-functor} holds because both $1_{W_0}^{\nat}(y_0,y_0)$ and $1_Z^{\nat}(z_0,z_0)$ are \emph{equal} to $1_{|y_0|}$, which may not be true when $\CQ$ is not integral. A simple counterexample in which Theorem \ref{hyperconvex=injective} fails is provided below.

Let $\CQ$ be a one-object quantaloid with three morphisms $\bot<1<\top$, where $1$ is the identity morphism (so that $\CQ$ is non-integral and involutive). It is easy to see that the one-object $\CQ$-category $\{1\}$ with hom $1$ is hyperconvex (note that $\top\circ\top=\top$). However, $\{1\}$ is not injective in $\QSymCat$. Indeed, let $\{\top\}$ be the one-object $\CQ$-category with hom $\top$, and let $\varnothing$ be the empty $\CQ$-category. Assuming $\{1\}$ is injective, there must be a $\CQ$-functor $h:\{\top\}\to\{1\}$, which does not exist (because $\top>1$ means that the unique map from $\{\top\}$ to $\{1\}$ is not $\CQ$-functorial).
$$\bfig
\qtriangle/->`->`-->/<600,400>[\varnothing`\{\top\}`\{1\};g`f`h]
\efig$$
\end{rem}

\section{Injective hulls of symmetric quantaloid-enriched categories} \label{Inj-hull-Sym-Q-Cat}

In this section we construct the injective hull (with respect to fully faithful $\CQ$-functors) of each symmetric $\CQ$-category in $\QSymCat$. To this end, for each symmetric $\CQ$-category $X$ we define
\begin{equation} \label{LX-TX-def}
\LX:=\{\mu\in\PX\mid\mu^{\circ}\circ\mu\leq 1_X^{\nat}\}\quad\text{and}\quad\TX:=\{\mu\in\PX\mid\mu^{\circ}=1_X^{\nat}\lda\mu\}
\end{equation}
as $\CQ$-subcategories of $\PX$, where $\TX$ is called the \emph{tight span} of $X$. The aim of this section is to show that $\TX$ is precisely the injective hull of $X$ in $\QSymCat$.

\begin{rem} \label{TX-def-QRel}
In order to facilitate the definition of the tight span of a partial metric space in Section \ref{Example}, we point out that a $\CQ$-relation $\mu:X\rto\{q\}$ satisfying 
$$\mu^{\circ}=1_X^{\nat}\lda\mu,\quad\text{or equivalently},\quad\mu=\mu^{\circ}\rda 1_X^{\nat},$$
must be a presheaf on $X$. Indeed, the above equality implies that
$$\mu\circ 1_X^{\nat}=(\mu^{\circ}\rda 1_X^{\nat})\circ 1_X^{\nat}\leq\mu^{\circ}\rda(1_X^{\nat}\circ 1_X^{\nat})=\mu^{\circ}\rda 1_X^{\nat}=\mu,$$
which shows that $\mu:X\oto\{q\}$ is a $\CQ$-distributor.
\end{rem}

\begin{lem} \label{TX-symmetric}
$\TX$ is a symmetric $\CQ$-category.
\end{lem}

\begin{proof}
Combining \eqref{lda-rda-circ} and \eqref{LX-TX-def}, it is straightforward to verify that
\begin{align*}
1_{\TX}^{\nat}(\lam,\mu)^{\circ}&=(\mu\lda\lam)^{\circ}=\lam^{\circ}\rda\mu^{\circ}=\lam^{\circ}\rda(1_X^{\nat}\lda\mu)\\
&=(\lam^{\circ}\rda 1_X^{\nat})\lda\mu=\lam\lda\mu=1_{\TX}^{\nat}(\mu,\lam)
\end{align*}
for all $\mu,\lam\in\TX$.
\end{proof}

Note that $\TX$ is non-empty as long as $X$ is non-empty. Indeed, for each symmetric $\CQ$-category $X$, it is easy to see that $1_X^{\nat}(-,x)\in\TX$ for all $x\in X$. In fact, there is a fully faithful $\CQ$-functor
\begin{equation} \label{y-def}
\sy_X:X\to\TX,\quad x\mapsto 1_X^{\nat}(-,x)
\end{equation}
that embeds $X$ into $\TX$, which is actually the codomain restriction of the \emph{Yoneda embedding} (cf. \cite{Stubbe2005,Shen2013a}). Recall that the \emph{Yoneda lemma} \cite{Stubbe2005} states that
\begin{equation} \label{Yoneda-lemma}
1_{\PX}^{\nat}(1_X^{\nat}(-,x),\mu)=\mu(x)
\end{equation}
for all $\mu\in\PX$ and $x\in X$. In particular, by restricting \eqref{Yoneda-lemma} to $\mu\in\TX$ we obtain that
\begin{equation} \label{Yoneda-lemma-TX}
(\sy_X)_{\nat}(-,\mu)=\mu.
\end{equation}

It is clear that $\TX\subseteq\LX$. We point out that $\TX$ consist of maximal elements in $\LX$; that is, if $\lam\in\TX$, $\mu\in\LX$ and $\lam\leq\mu$, then
$$\mu^{\circ}\leq 1_X^{\nat}\lda\mu\leq 1_X^{\nat}\lda\lam=\lam^{\circ},$$
which forces $\mu=\lam$. Conversely:

\begin{lem} \label{tmu-maximum}
For each $\mu\in\LX$, there exists $\tmu\in\TX$ such that $\mu\leq\tmu$.
\end{lem}

\begin{proof}
Let us consider the set
$$Y:=\{\lam\in\LX\mid\mu\leq\lam\}.$$
Since $\mu\in Y$, it is clear that $Y\neq\varnothing$. Note that every chain $\{\lam_j\mid j\in J\}\subseteq Y$ has an upper bound in $Y$, given by $\bv\limits_{j\in J}\lam_j$. By Zorn's lemma, $Y$ has a maximal element $\tmu$. We claim that $\tmu\in\TX$.

We proceed by contradiction. Assume that $\tmu\not\in\TX$. Then there exists $z\in X$ such that
$$\tmu^{\circ}(z)<1_X^{\nat}(-,z)\lda\tmu,$$
which is equivalent to
$$\tmu(z)<\tmu^{\circ}\rda(1_X^{\nat}(-,z))^{\circ}=\tmu^{\circ}\rda 1_X^{\nat}(z,-).$$
Define a presheaf $\mu_0$ on $X$ given by
$$\mu_0:=\tmu\vee((\tmu^{\circ}\rda 1_X^{\nat}(z,-))\circ 1_X^{\nat}(-,z)).$$
Then $\mu_0\in\LX$. Indeed, since
$$\mu_0^{\circ}=\tmu^{\circ}\vee(1_X^{\nat}(z,-)\circ(1_X^{\nat}(-,z)\lda\tmu)),$$
the inequality $\mu_0^{\circ}\circ\mu_0\leq 1_X^{\nat}$ is a direct consequence of
\begin{align*}
&\tmu^{\circ}\circ\tmu\leq 1_X^{\nat},\\
&\tmu^{\circ}\circ(\tmu^{\circ}\rda 1_X^{\nat}(z,-))\circ 1_X^{\nat}(-,z)\leq 1_X^{\nat}(z,-)\circ 1_X^{\nat}(-,z)\leq 1_X^{\nat},\\
&1_X^{\nat}(z,-)\circ(1_X^{\nat}(-,z)\lda\tmu)\circ\tmu\leq 1_X^{\nat}(z,-)\circ 1_X^{\nat}(-,z)\leq 1_X^{\nat}
\end{align*}
and
\begin{align*}
&1_X^{\nat}(z,-)\circ(1_X^{\nat}(-,z)\lda\tmu)\circ(\tmu^{\circ}\rda 1_X^{\nat}(z,-))\circ 1_X^{\nat}(-,z)\\
\leq{}&1_X^{\nat}(z,-)\circ 1_{|z|}\circ 1_X^{\nat}(-,z)&(\CQ\ \text{is integral})\\
={}&1_X^{\nat}(z,-)\circ 1_X^{\nat}(-,z)\leq 1_X^{\nat}.
\end{align*}
But $\tmu<\mu_0$, because
$$\tmu(z)<\tmu^{\circ}\rda 1_X^{\nat}(z,-)=(\tmu^{\circ}\rda 1_X^{\nat}(z,-))\circ 1_X^{\nat}(z,z)\leq\mu_0(z),$$
contradicting to the maximality of $\tmu$.
\end{proof}

With Lemma \ref{tmu-maximum} we are able to prove the hyperconvexity of $\TX$:

\begin{lem} \label{TX-hyperconvex}
$\TX$ is hyperconvex.
\end{lem}

\begin{proof}
Let $\Theta$ be a presheaf on $\TX$ satisfying $\Theta^{\circ}\circ\Theta\leq 1_{\TX}^{\nat}$. Let $\lam\in\PX$ be given by
$$\lam=\bv_{\mu\in\TX}\Theta(\mu)\circ\mu.$$
Then $\lam\in\LX$, since
\begin{align*}
\lam^{\circ}\circ\lam&=\Big(\bv_{\mu\in\TX}\Theta(\mu)\circ\mu\Big)^{\circ}\circ\Big(\bv_{\mu'\in\TX}\Theta(\mu')\circ\mu'\Big)\\
&=\bv_{\mu\in\TX}\bv_{\mu'\in\TX}\mu^{\circ}\circ\Theta^{\circ}(\mu)\circ\Theta(\mu')\circ\mu'\\
&\leq\bv_{\mu\in\TX}\bv_{\mu'\in\TX}\mu^{\circ}\circ 1_{\TX}^{\nat}(\mu',\mu)\circ\mu'&(\Theta^{\circ}\circ\Theta\leq 1_{\TX}^{\nat})\\
&=\bv_{\mu\in\TX}\bv_{\mu'\in\TX}\mu^{\circ}\circ(\mu\lda\mu')\circ\mu'\\
&\leq\bv_{\mu\in\TX}\mu^{\circ}\circ\mu\\
&\leq 1_X^{\nat}.&(\TX\subseteq\LX\ \text{and Equations \eqref{LX-TX-def}})
\end{align*}
By Lemma \ref{tmu-maximum}, there exists $\tlam\in\TX$ such that $\lam\leq\tlam$. It follows that
$$\Theta(\mu)\leq\lam\lda\mu\leq\tlam\lda\mu=1_{\TX}^{\nat}(\mu,\tlam)$$
for all $\mu\in\TX$, which completes the proof.
\end{proof}

A fully faithful $\CQ$-functor $f:X\to Y$ between symmetric $\CQ$-categories is \emph{essential} in $\QSymCat$ if, for each symmetric $\CQ$-category $Z$, a $\CQ$-functor $g:Y\to Z$ is fully faithful whenever $gf:X\to Z$ is fully faithful.

A symmetric $\CQ$-category $Y$ is the \emph{injective hull} of $X$ in $\QSymCat$, if $Y$ is injective in $\QSymCat$ and there exists an essential fully faithful $\CQ$-functor $f:X\to Y$. It is well known that injective hulls are \emph{essentially unique}; that is, if $Y'$ is another injective hull of $X$ in $\QSymCat$ (with an essential fully faithful $\CQ$-functor $f':X\to Y'$), then there exists an isomorphism $g:Y\to Y'$ in $\QSymCat$ with $f'=gf$.

For each symmetric $\CQ$-category $X$, since $\TX$ is injective in $\QSymCat$ by Lemma \ref{TX-hyperconvex} and Theorem \ref{hyperconvex=injective}, it is actually the injective hull of $X$ in $\QSymCat$:

\begin{thm} \label{injective-hull}
Let $X$ be a symmetric $\CQ$-category. Then the fully faithful $\CQ$-functor $\sy_X:X\to\TX$ is essential. Hence, the tight span $\TX$ is the injective hull of $X$ in $\QSymCat$.
\end{thm}

\begin{proof}
Let $Y$ be a symmetric $\CQ$-category. Let $g:\TX\to Y$ be a $\CQ$-functor such that $g\sy_X:X\to Y$ is fully faithful. Note that for any $\mu,\lam\in\TX$, from \eqref{LX-TX-def} and \eqref{Yoneda-lemma-TX} we see that
$$\sy_X^{\nat}(\mu,-)=(\sy_X)_{\nat}(-,\mu)^{\circ}=\mu^{\circ}=1_X^{\nat}\lda\mu=1_X^{\nat}\lda(\sy_X)_{\nat}(-,\mu)$$
and
$$1_{\TX}^{\nat}(\mu,\lam)=\lam\lda\mu=(\sy_X)_{\nat}(-,\lam)\lda(\sy_X)_{\nat}(-,\mu);$$
that is, $\sy_X^{\nat}=1_X^{\nat}\lda(\sy_X)_{\nat}$ and $1_{\TX}^{\nat}=(\sy_X)_{\nat}\lda(\sy_X)_{\nat}$. By Lemmas \ref{fully-faithful-graph} and \ref{graph-cograph-calc}\ref{graph-cograph-calc:composition}, the fully faithfulness of $g\sy_X$ means that
$$\sy_X^{\nat}\circ g^{\nat}\circ g_{\nat}\circ(\sy_X)_{\nat}=(g\sy_X)^{\nat}\circ(g\sy_X)_{\nat}=1_X^{\nat},$$
and consequently
$$g^{\nat}\circ g_{\nat}\leq\sy_X^{\nat}\rda(1_X^{\nat}\lda(\sy_X)_{\nat})=\sy_X^{\nat}\rda \sy_X^{\nat}=((\sy_X)_{\nat}\lda(\sy_X)_{\nat})^{\circ}=(1_{\TX}^{\nat})^{\circ}=1_{\TX}^{\nat}.$$
As the reverse inequality is an immediate consequence of $g_{\nat}\dv g^{\nat}$, we deduce that $g^{\nat}\circ g_{\nat}=1_{\TX}^{\nat}$, and therefore the fully faithfulness of $g$ follows from Lemma \ref{fully-faithful-graph}.
\end{proof}

\section{Essential embeddings of symmetric quantaloid-enriched categories} \label{Essential-embed-Sym-Q-Cat}

Recall that a $\CQ$-functor $f:X\to Y$ is \emph{dense} (resp. \emph{codense}) (cf. \cite[Proposition 4.12]{Lai2017}) if
\begin{equation} \label{dense-def}
1_Y^{\nat}=f_{\nat}\lda f_{\nat}\quad(\text{resp.}\ 1_Y^{\nat}=f^{\nat}\rda f^{\nat}).
\end{equation}
Note that the notions of density and codensity coincide when $Y$ is symmetric: by Equations \eqref{f-graph-circ},
\begin{equation} \label{dense=codense}
1_Y^{\nat}=f_{\nat}\lda f_{\nat}\iff (1_Y^{\nat})^{\circ}=(f_{\nat}\lda f_{\nat})^{\circ}\iff 1_Y^{\nat}=f^{\nat}\rda f^{\nat}.
\end{equation}
The aim of this section is to characterize essential fully faithful $\CQ$-functors in $\QSymCat$ (i.e., ``essential embeddings'' of symmetric $\CQ$-categories) through their (co)density. 

\begin{lem} \label{dense-TX}
Let $f:X\to Y$ be a $\CQ$-functor between symmetric $\CQ$-categories. If $f$ is dense and fully faithful, then $f_{\nat}(-,y)\in\TX$ for all $y\in Y$.
\end{lem}

\begin{proof}
Combining Equation \eqref{dense-def}, Lemmas  \ref{graph-cograph-calc}\ref{graph-cograph-calc:f-nat-bracket} and \ref{fully-faithful-graph} we compute that
\begin{equation} \label{f-nat-y-TX}
f^{\nat}=f^{\nat}\circ(f_{\nat}\lda f_{\nat})=(f^{\nat}\circ f_{\nat})\lda f_{\nat}=1_X^{\nat}\lda f_{\nat};
\end{equation}
that is,
$$f_{\nat}(-,y)^{\circ}=f^{\nat}(y,-)=1_X^{\nat}\lda f_{\nat}(-,y),$$
i.e., $f_{\nat}(-,y)\in\TX$, for all $y\in Y$.
\end{proof}

The following Lemma \ref{TX-LX-retract} strengthens Lemma \ref{tmu-maximum} by revealing that $\TX$ is a retract of the symmetrization of $\LX$ (see \eqref{symmetrization-def}):

\begin{lem} \label{TX-LX-retract}
For each symmetric $\CQ$-category $X$, there exists a $\CQ$-functor 
$$m_X:(\LX)_{\sfs}\to\TX$$ 
whose restriction on $\TX$ is $1_{\TX}$, such that $\mu\leq m_X\mu$ for all $\mu\in\LX$.
\end{lem}

\begin{proof}
Since $\TX$ is symmetric, $(\LX)_{\sfs}$ is a $\CQ$-supercategory of $\TX$. Since $\TX$ is an absolute retract in $\QSymCat$ by the hyperconvexity of $\TX$ (Lemma \ref{TX-hyperconvex}) and Theorem \ref{hyperconvex=injective}, there exists a $\CQ$-functor $m_X:(\LX)_{\sfs}\to\TX$ with $m_X|_{\TX}=1_{\TX}$.

To show that $\mu\leq m_X\mu$ for any $\mu\in\LX$, note that $\mu\leq\mu^{\circ}\rda 1_X^{\nat}$ by Equations \eqref{LX-TX-def}, and consequently,
\begin{align*}
\mu(x)&=\mu(x)\wedge(\mu^{\circ}\rda 1_X^{\nat}(x,-))&(\mu\leq\mu^{\circ}\rda 1_X^{\nat})\\
&=1_{\PX}^{\nat}(1_X^{\nat}(-,x),\mu)\wedge 1_{\PX}^{\nat}(\mu,1_X^{\nat}(-,x))^{\circ}&(\text{Equation \eqref{Yoneda-lemma}})\\
&=1_{(\LX)_{\sfs}}^{\nat}(1_X^{\nat}(-,x),\mu)\\
&\leq 1_{\TX}^{\nat}(m_X(1_X^{\nat}(-,x)),m_X\mu)&(m_X\ \text{is a}\ \CQ\text{-functor})\\
&=1_{\TX}^{\nat}(1_X^{\nat}(-,x),m_X\mu)&(1_X^{\nat}(-,x)\in\TX\ \text{and}\ m_X|_{\TX}=1_{\TX})\\
&=(m_X\mu)(x)&(\text{Equation \eqref{Yoneda-lemma}})
\end{align*}
for all $x\in X$.
\end{proof}

\begin{lem} \label{TX-TY-iso}
Let $f:X\to Y$ be a fully faithful $\CQ$-functor between symmetric $\CQ$-categories. If $f$ is essential, then $\TX$ and $\TY$ are isomorphic $\CQ$-categories, with an isomorphism in $\QSymCat$ given by
$$(-)\circ f_{\nat}:\TY\to\TX,\quad \lam\mapsto\lam\circ f_{\nat}.$$
\end{lem}

\begin{proof}
Since $f$ is fully faithful, for each $\lam\in\TY$, it follows from \eqref{LX-TX-def} and Lemma \ref{fully-faithful-graph} that
$$(\lam\circ f_{\nat})^{\circ}\circ(\lam\circ f_{\nat})=f^{\nat}\circ\lam^{\circ}\circ\lam\circ f_{\nat}\leq f^{\nat}\circ f_{\nat}=1_X^{\nat};$$
that is, $\lam\circ f_{\nat}\in\LX$. Thus, it is easy to see that $(-)\circ f_{\nat}$ defines a $\CQ$-functor from $\TY$ to $\LX$, and consequently a $\CQ$-functor from $\TY$ to $(\LX)_{\sfs}$ (see Lemma \ref{QSymCat-coref-QCat}), which induces a $\CQ$-functor
$$g:=\big(\TY\to^{(-)\circ f_{\nat}}(\LX)_{\sfs}\to^{m_X}\TX\big).$$
Conversely, note that for each $\mu\in\TX$, 
$$(\mu\circ f^{\nat})^{\circ}\circ(\mu\circ f^{\nat})=f_{\nat}\circ\mu^{\circ}\circ\mu\circ f^{\nat}\leq f_{\nat}\circ f^{\nat}\leq 1_Y^{\nat};$$
that is, $\mu\circ f^{\nat}\in\LY$. Thus, similarly, there is a $\CQ$-functor
$$h:=\big(\TX\to^{(-)\circ f^{\nat}}(\LY)_{\sfs}\to^{m_Y}\TY\big).$$

First, $g$ is fully faithful. Since $\sy_Y:Y\to\TY$ is essential (Theorem \ref{injective-hull}), it is easy to see that the composite $\sy_Y f:X\to\TY$ is essential. Since
$$(\sy_Y fx)\circ f_{\nat}=1_Y^{\nat}(-,fx)\circ f_{\nat}=1_Y^{\nat}(f-,fx)=1_X^{\nat}(-,x)$$
by Lemma \ref{graph-cograph-calc}\ref{graph-cograph-calc:phi-f-g} and the fully faithfulness of $f$, we obtain that
$$g\sy_Y fx=m_X((\sy_Y fx)\circ f_{\nat})=m_X(1_X^{\nat}(-,x))=1_X^{\nat}(-,x)$$
for all $x\in X$. Consequently, 
$$1_X^{\nat}(x,x')=1_X^{\nat}(-,x')\lda 1_X^{\nat}(-,x)=1_{\TX}^{\nat}(g\sy_Y fx,g\sy_Y fx')$$
for all $x,x'\in X$, showing that $g\sy_Y f$ is fully faithful. Hence, the fully faithfulness of $g$ follows from the essentiality of $\sy_Y f$.

Second, $g$ is surjective. To this end, we show that $gh\mu=\mu$ for all $\mu\in\TX$. Indeed, since $\mu\circ f^{\nat}\in\LY$, it holds that
$$\mu=\mu\circ f^{\nat}\circ f_{\nat}\leq m_Y(\mu\circ f^{\nat})\circ f_{\nat}$$
by the fully faithfulness of $f$ (Lemma \ref{fully-faithful-graph}) and Lemma \ref{TX-LX-retract}. Thus, the maximality of presheaves in $\TX$ forces 
\begin{equation} \label{mu-leq-mY-mu-f-nat}
\mu=m_Y(\mu\circ f^{\nat})\circ f_{\nat},
\end{equation}
and consequently
$$gh\mu=m_X(m_Y(\mu\circ f^{\nat})\circ f_{\nat})=m_X\mu=\mu.$$

Therefore, as a fully faithful and surjective $\CQ$-functor between separated $\CQ$-categories, $g:\TY\to\TX$ is necessarily an isomorphism in $\QSymCat$. Since we have already proved that $gh=1_{\TX}$, $h:\TX\to\TY$ must be the inverse of $g$; that is, $hg=1_{\TY}$. Thus, for every $\lam\in\TY$, replacing $\mu$ with $g\lam$ in \eqref{mu-leq-mY-mu-f-nat} (note that \eqref{mu-leq-mY-mu-f-nat} holds for all $\mu\in\TX$) we deduce that $$g\lam=m_Y((g\lam)\circ f^{\nat})\circ f_{\nat}=(hg\lam)\circ f_{\nat}=\lam\circ f_{\nat},$$ 
which completes the proof.
\end{proof}

Now we are ready to present the main result of this section:

\begin{thm} \label{f-essential}
Let $f:X\to Y$ be a fully faithful $\CQ$-functor between symmetric $\CQ$-categories. Then the following statements are equivalent:
\begin{enumerate}[label={\rm(\roman*)}]
\item \label{f-essential:essential} $f$ is essential in $\QSymCat$.
\item \label{f-essential:dense} $f$ is dense.
\item \label{f-essential:codense} $f$ is codense.
\end{enumerate}
\end{thm}

\begin{proof}
As \ref{f-essential:dense}$\iff$\ref{f-essential:codense} automatically holds by \eqref{dense=codense}, it remains to prove \ref{f-essential:essential}$\implies$\ref{f-essential:dense} and \ref{f-essential:dense}+\ref{f-essential:codense}$\implies$\ref{f-essential:essential}.

\ref{f-essential:essential}$\implies$\ref{f-essential:dense}: Suppose that $f$ is essential. By Lemma \ref{TX-TY-iso},
$$g:=\big(Y\to^{\sy_Y}\TY\to^{(-)\circ f_{\nat}}\TX\big)$$
is a well-defined $\CQ$-functor. Note that 
$$gfx=1_Y^{\nat}(-,fx)\circ f_{\nat}=1_Y^{\nat}(f-,fx)=1_X^{\nat}(-,x)$$
for all $x\in X$, the fully faithfulness of $gf$ follows immediately. Hence, $g$ is also fully faithful, which means that
$$1_Y^{\nat}(y,y')=1_{\TX}^{\nat}((\sy_Y y)\circ f_{\nat},(\sy_Y y')\circ f_{\nat})=1_{\TX}^{\nat}(f_{\nat}(-,y),f_{\nat}(-,y'))=f_{\nat}(-,y')\lda f_{\nat}(-,y)$$
for all $y,y'\in Y$; that is, $1_Y^{\nat}=f_{\nat}\lda f_{\nat}$.

\ref{f-essential:dense}+\ref{f-essential:codense}$\implies$\ref{f-essential:essential}: Suppose that $Z$ is symmetric, $g:Y\to Z$ is a $\CQ$-functor and $gf:X\to Z$ is fully faithful. Applying Lemma \ref{fully-faithful-graph} to the fully faithful $\CQ$-functor $gf$, we have 
$$f^{\nat}\circ g^{\nat}\circ g_{\nat}\circ f_{\nat}=(gf)^{\nat}\circ(gf)_{\nat}=1_X^{\nat}.$$
Thus, with Equation \eqref{f-nat-y-TX} in Lemma \ref{dense-TX} and the definition of codensity (Equation \eqref{dense-def}) it is easy to compute that
$$g^{\nat}\circ g_{\nat}\leq f^{\nat}\rda(1_X^{\nat}\lda f_{\nat})=f^{\nat}\rda f^{\nat}=1_Y^{\nat}.$$
As the reverse inequality is an immediate consequence of $g_{\nat}\dv g^{\nat}$, we deduce that $g^{\nat}\circ g_{\nat}=1_Y^{\nat}$, and therefore the fully faithfulness of $g$ follows from Lemma \ref{fully-faithful-graph}.
\end{proof}

\begin{rem}
Lemma \ref{TX-TY-iso} is pivotal in the proof of ``\ref{f-essential:essential}$\implies$\ref{f-essential:dense}'' of Theorem \ref{f-essential}, where deriving that $f_{\nat}(-,y)\in\TX$ for all $y\in Y$ is the crucial step. The proof of Lemma \ref{TX-TY-iso} is highly non-trivial, where we invoke some ideas from \cite[Subsection 1.13]{Dress1984}. However, our proof of Lemma \ref{TX-LX-retract}, which provides the important retraction needed in the proof of Lemma \ref{TX-TY-iso},  deviates from \cite{Dress1984}: by proving the hyperconvexity of $\TX$ (Lemma \ref{TX-hyperconvex}) prior to Lemma \ref{TX-LX-retract}, we are able to immediately deduce the existence of $m_X:(\LX)_{\sfs}\to\TX$.
\end{rem}

\begin{rem} \label{f-essential-comparison}
In the case that $\CQ=[0,\infty]$ is the Lawvere quantale (see Example \ref{quantale-exmp}\ref{quantale-exmp:Lawvere}), if we compare Theorem \ref{f-essential} with \cite[Example 9.13(7)]{Adamek1990} or \cite[Theorem 2]{Herrlich1992a}, it seems that our theorem should be parallelly stated as: A fully faithful $\CQ$-functor $f:X\to Y$ between symmetric $\CQ$-categories is essential if, and only if, $f$ satisfies the following two conditions:
\begin{enumerate}[label={\rm(\alph*)}]
\item \label{f-essential-Met:dense} $1_Y^{\nat}=f_{\nat}\lda f_{\nat}$, i.e., $f$ is dense;
\item \label{f-essential-Met:TX} $f^{\nat}=1_X^{\nat}\lda f_{\nat}$, i.e., $f_{\nat}(-,y)\in\TX$ for all $y\in Y$.
\end{enumerate}
In fact, our Lemma \ref{dense-TX} ensures that the condition \ref{f-essential-Met:TX} is contained in \ref{f-essential-Met:dense}, and therefore the (co)density of a fully faithful $\CQ$-functor is sufficient to derive its essentiality. 
\end{rem}

\section{Example: injective partial metric spaces} \label{Example}

A \emph{(unital) quantale} \cite{Mulvey1986,Rosenthal1990} is a one-object quantaloid. In this section, we let
$$\sQ=(\sQ,\otimes,1,{}^{\circ})$$
denote an \emph{integral} and \emph{involutive} quantale. Explicitly: 
\begin{itemize}
\item $(\sQ,\otimes,1)$ is a monoid;
\item $\sQ$ is equipped with the structure of a complete lattice (with the unit $1$ being the top element);
\item the multiplication $\otimes$ preserves arbitrary suprema on both sides;
\item $(-)^{\circ}:\sQ\to\sQ$ is a sup-preserving map with 
$$q^{\circ\circ}=q\quad\text{and}\quad (p\otimes q)^{\circ}=q^{\circ}\otimes p^{\circ}$$
for all $p,q\in\sQ$. 
\end{itemize}
To avoid confusion, we denote the left and right implications in $\sQ$ by $\ldd$ and $\rdd$, respectively, which satisfy
$$p\otimes q\leq r\iff p\leq r\ldd q\iff q\leq p\rdd r$$
for all $p,q,r\in\sQ$. We say that:
\begin{itemize}
\item $\sQ$ is \emph{commutative}, if $p\otimes q=q\otimes p$ for all $p,q\in\sQ$, in which case we write
$$p\ra q:=q\ldd p=p\rdd q$$
for all $p,q\in\sQ$. Note that every commutative quantale $\sQ$ is involutive, with a trivial involution given by the identity map $1_{\sQ}:\sQ\to\sQ$.
\item $\sQ$ is \emph{divisible}, if
\begin{equation} \label{Q-divisible-def}
(u\ldd q)\otimes q=u=q\otimes(q\rdd u)
\end{equation}
whenever $u\leq q$ in $\sQ$. Note that every quantale satisfying \eqref{Q-divisible-def} is necessarily integral (cf. \cite[Proposition 2.1]{Pu2012} and \cite[Proposition 3.1]{Tao2014}).
\end{itemize}

\begin{exmp} \label{quantale-exmp}
We list here some examples of integral and involutive quantales:
\begin{enumerate}[label=(\arabic*)]
\item \label{quantale-exmp:Lawvere} The Lawvere quantale $[0,\infty]=([0,\infty],+,0)$ \cite{Lawvere1973} is commutative and divisible, where $[0,\infty]$ is the extended non-negative real line equipped with the order ``$\geq$'', and ``$+$'' is the usual addition extended via
$$p+\infty=\infty+p=\infty$$
to $[0,\infty]$. The implication in $[0,\infty]$ is given by\footnote{In order to eliminate ambiguity, we make the convention that all the symbols (e.g., $\leq$, $\vee$, $\max$, $\sup$, etc.) connecting (extended) real numbers always refer to the standard order on $[0,\infty]$, although the quantale $[0,\infty]$ is equipped with the reverse order ``$\geq$'' of (extended) real numbers.}
\begin{equation} \label{0-infty-imp}
p\ra q=\max\{0,q-p\}
\end{equation}
for all $p,q\in[0,\infty]$. In this case, symmetric $[0,\infty]$-categories are given by maps $\al:X\times X\to[0,\infty]$ such that 
\begin{itemize}
\item $\al(x,x)=0$,
\item $\al(x,y)=\al(y,x)$,
\item $\al(x,z)\leq\al(x,y)+\al(y,z)$
\end{itemize}
for all $x,y,z\in X$, and they become (classical) \emph{metric spaces} if, moreover,
$$\al(x,y)<\infty\quad\text{and}\quad x=y\iff\al(x,y)=0$$
for all $x,y\in X$.
\item \label{quantale-exmp:frame} Every \emph{frame} $\Om=(\Om,\wedge,\top)$ is a commutative and divisible quantale.
\item \label{quantale-exmp:t-norm} Every \emph{complete BL-algebra} \cite{Hajek1998} is a commutative and divisible quantale. In particular, the unit interval $[0,1]$ equipped with a \emph{continuous t-norm} \cite{Klement2000} is a commutative and divisible quantale, and so is every \emph{complete MV-algebra} \cite{Chang1958}.
\item \label{quantale-exmp:nil-min} The unit interval $[0,1]$ equipped with the \emph{nilpotent minimum t-norm} \cite{Klement2000} is a commutative, integral and non-divisible quantale.
\item \label{quantale-exmp:Sup} Let $\Sup[0,1]$ denote the set of $\sup$-preserving maps on the unit interval $[0,1]$. Then 
$$\Sup[0,1]=(\Sup[0,1],\circ,1_{[0,1]})$$ 
is a non-commutative, non-integral and involutive quantale \cite{Eklund2018}, whose multiplication is given by the composition $\circ$ of maps, and the unit is given by the identity map $1_{[0,1]}$ on $[0,1]$. An involution on $\Sup[0,1]$ is given by
$$f^{\circ}:[0,1]\to[0,1],\quad f^{\circ}(x)=1-f^{\star}(1-x)$$
for all $f\in\Sup[0,1]$, where $f^{\star}:[0,1]\to[0,1]$ is the right adjoint of $f$. Furthermore, it is straightforward to verify that 
$$f\leq 1_{[0,1]}\iff f^{\circ}\leq 1_{[0,1]}$$ 
for all $f\in\Sup[0,1]$, and consequently, 
$$\Sup[0,1]_{\leq 1_{[0,1]}}=\{f\in\Sup[0,1]\mid f\leq 1_{[0,1]}\}$$
is a subquantale of $\Sup[0,1]$ that is non-commutative, non-divisible, integral and involutive.
\end{enumerate}
\end{exmp}

Since integral and involutive quantales are precisely integral and involutive quantaloids with only one object, our results in Sections \ref{Inj-Sym-Q-Cat}--\ref{Essential-embed-Sym-Q-Cat} generalize those in \cite{Jawhari1986,Kabil2020}. In particular, by setting $\CQ=[0,\infty]$, our Theorems \ref{hyperconvex=injective} and \ref{injective-hull} are reduced to Theorems \ref{Met-injective} and \ref{Met-injective-hull} (by the aid of Remarks \ref{hyperconvex-def-QRel} and \ref{TX-def-QRel}), respectively, which characterize injective objects and injective hulls in the category $\Met$ of (classical) metric spaces and non-expansive maps. Moreover, as we point out in Remark \ref{f-essential-comparison}, our Theorem \ref{f-essential} generalizes and refines the existing characterization of essential embeddings of metric spaces.

\begin{rem}
Since the implication in the quantale $[0,\infty]$ is given by \eqref{0-infty-imp}, while applying \eqref{LX-TX-def} and Remark \ref{TX-def-QRel} to the case of $\CQ=[0,\infty]$, in Theorem \ref{Met-injective-hull} it seems that $\sT(X,\al)$ should consist of maps $\mu:X\to[0,\infty]$ satisfying
$$\mu(x)=\max\{0,\sup\limits_{y\in X}(\al(x,y)-\mu(y))\}$$
for all $x\in X$, and
$$\si(\mu,\lam)=\max\{0,\sup\limits_{x\in X}(\mu(x)-\lam(x))\}$$
for all $\mu,\lam\in\sT(X,\al)$. Indeed, it is straightforward to verify that
$$\mu(x)=\sup\limits_{y\in X}(\al(x,y)-\mu(y))\iff\mu(x)=\max\{0,\sup\limits_{y\in X}(\al(x,y)-\mu(y))\}$$
for every map $\mu:X\to[0,\infty]$ and $x\in X$, and 
$$\sup\limits_{x\in X}(\mu(x)-\lam(x))=\sup\limits_{x\in X}(\lam(x)-\mu(x))\geq 0$$
for all $\mu,\lam\in\sT(X,\al)$. So, Theorem \ref{Met-injective-hull} is indeed a special case of our Theorem \ref{injective-hull}.
\end{rem}

In fact, the potential applications of our results for a general (small, integral and involutive) quantaloid are far beyond the one-object case. It is well known that $\sQ$ gives rise to a quantaloid $\DQ$ of \emph{diagonals} of $\sQ$ \cite{Hoehle2011a,Pu2012,Stubbe2014}, which consists of the following data:
\begin{itemize}
\item objects of $\DQ$ are elements $p,q,r,\dots$ of $\sQ$;
\item for $p,q\in\sQ$, a morphism $u:p\to q$ in $\DQ$, called a \emph{diagonal} from $p$ to $q$, is an element $u\in\sQ$ with
    \begin{equation} \label{diagonal-def}
    (u\ldd p)\otimes p=u=q\otimes(q\rdd u);
    \end{equation}
\item for diagonals $u:p\to q$, $v:q\to r$, the composition $v\circ u:p\to r$ is given by
    \begin{equation} \label{diagonal-comp}
    v\circ u=(v\ldd q)\otimes q\otimes(q\rdd u)=(v\ldd q)\otimes u=v\otimes(q\rdd u);
    \end{equation}
\item $q:q\to q$ is the identity diagonal on $q$.
\end{itemize}

Since $\sQ$ is integral, Equation \eqref{diagonal-def} necessarily forces $u\leq p\wedge q$ for all morphisms $u:p\to q$ in $\DQ$, which guarantees the integrality of $\DQ$. In particular:
\begin{itemize}
\item If $\sQ$ is a divisible quantale, then $u$ satisfies \eqref{diagonal-def} if, and only if, $u\leq p\wedge q$.
\item If $\sQ=\Om$ is a frame, then the composition \eqref{diagonal-comp} becomes
\begin{equation} \label{DOm-comp}
v\circ u=v\wedge u
\end{equation}
for all diagonals $u:p\to q$, $v:q\to r$ \cite{Walters1981}.
\end{itemize}

Note also that $u:p\to q$ is a diagonal if, and only if, $u^{\circ}:q^{\circ}\to p^{\circ}$ is a diagonal. Hence, the full subquantaloid $\DQ^{\circ}$ of $\DQ$ given by
$$\ob\DQ^{\circ}:=\{q\in\sQ\mid q^{\circ}=q\}$$
is a small, integral and involutive quantaloid, with the involution lifted from $\sQ$. In particular, we have
$$\DQ^{\circ}=\DQ$$
if $\sQ$ is commutative.

\begin{exmp} \label{Q-set}
A symmetric $\DQ^{\circ}$-category is exactly a \emph{$\sQ$-set} in the sense of H{\"o}hle--Kubiak (see \cite[Definition 2.1]{Hoehle2011a}), which may be described as a set $X$ equipped with a map $\al:X\times X\to\sQ$ such that
such that
\begin{enumerate}[label=(S\arabic*)]
\item \label{sim-def:str} $\al(x,y)\leq\al(x,x)\wedge\al(y,y)$,
\item \label{sim-def:sym} $\al(x,y)=\al(y,x)^{\circ}$,
\item \label{sim-def:div} $\al(x,y)=(\al(x,y)\ldd\al(x,x))\otimes\al(x,x)$,
\item \label{sim-def:tran} $(\al(y,z)\ldd\al(y,y))\otimes\al(x,y)\leq\al(x,z)$
\end{enumerate}
for all $x,y,z\in X$, where the type map $|\text{-}|:X\to\sQ$ is given by 
$$|x|=\al(x,x)$$ 
for all $x\in X$; for details we refer to \cite[Proposition 6.3]{Hoehle2011a} and \cite[Theorem 4.5]{Lai2020}. Note that \ref{sim-def:str} is actually subsumed by \ref{sim-def:div} since our $\sQ$ is integral, and \ref{sim-def:str} is equivalent to \ref{sim-def:div} when $\sQ$ is divisible. In particular: 
\begin{enumerate}[label=(\arabic*)]
\item \label{Q-set:divisible} If $\sQ$ is a commutative and divisible quantale, then a $\sQ$-set $(X,\al)$ is given by a map $\al:X\times X\to\sQ$ such that
\begin{itemize}
\item $\al(x,y)\leq\al(x,x)\wedge\al(y,y)$,
\item $\al(x,y)=\al(y,x)$,
\item $\al(y,z)\otimes(\al(y,y)\ra\al(x,y))\leq\al(x,z)$
\end{itemize}
for all $x,y,z\in X$.
\item \label{Q-set:metric} If $[0,\infty]$ is the Lawvere quantale, then a $[0,\infty]$-set $(X,\al)$ is given by a map $\al:X\times X\to[0,\infty]$ such that
\begin{itemize}
\item $\al(x,x)\vee\al(y,y)\leq\al(x,y)$,
\item $\al(x,y)=\al(y,x)$,
\item $\al(x,z)\leq\al(x,y)-\al(y,y)+\al(y,z)$
\end{itemize}
for all $x,y,z\in X$; that is, a (slightly generalized) \emph{partial metric space} \cite{Matthews1994,Bukatin2009,Hoehle2011a,Pu2012,Stubbe2014,Hofmann2018a}. We remind the readers that the notion of ``partial metric'' originally introduced by Matthews \cite{Matthews1994} additionally requires $\al$ to satisfy
$$\al(x,y)<\infty\quad\text{and}\quad x=y\iff\al(x,x)=\al(y,y)=\al(x,y)$$
for all $x,y\in X$.
\item \label{Q-set:frame} If $\Om$ is a frame, then an $\Om$-set \cite{Fourman1979,Borceux1994c} $(X,\al)$ is given by a map $\al:X\times X\to\Om$ such that
\begin{itemize}
\item $\al(x,y)=\al(y,x)$,
\item $\al(y,z)\wedge\al(x,y)\leq\al(x,z)$
\end{itemize}
for all $x,y,z\in X$.
\end{enumerate}
\end{exmp}

Therefore, for every integral and involutive quantale $\sQ$, our results of injective objects, injective hulls and essential embeddings in $\QSymCat$ can be seamlessly applied to the category of $\sQ$-sets, giving rise to a theory of \emph{injective $\sQ$-sets}. In the rest of this section, we derive \emph{injective partial metric spaces} as an example, where we consider the category
$$\ParMet:=\sD[0,\infty]\text{-}{\bf SymCat}$$
of (slightly generalized) partial metric spaces (as defined in Example \ref{Q-set}\ref{Q-set:metric}) and \emph{non-expansive maps} (i.e., $\sD[0,\infty]$-functors). 

Let $(X,\al)$ be a partial metric space. By Remark \ref{hyperconvex-def-QRel}, Definition \ref{hyperconvex-def} may be formulated in the context of partial metric spaces as:

\begin{defn} \label{parmet-hyperconvex-def}
A partial metric space $(X,\al)$ is \emph{hyperconvex} if, for every map $\mu:X\to[r,\infty]$  $(r\in[0,\infty])$ satisfying
$$\al(x,x)\leq\mu(x)\quad\text{and}\quad\al(x,y)\leq\mu(x)-r+\mu(y)$$
for all $x,y\in X$, there exists $z\in X$ such that
$$\al(x,z)\leq\mu(x)$$
for all $x\in X$.
\end{defn}

\begin{rem}
Let us restate Definition \ref{parmet-hyperconvex-def} in a geometric way. A partial metric space $(X,\al)$ is hyperconvex if, for every $r\in[0,\infty]$ and every family $\{(x_j,r_j)\}_{j\in J}$ of pairs with $x_j\in X$ and $r_j\in[r,\infty]$ satisfying
$$\al(x_j,x_j)\leq r_j\quad\text{and}\quad\al(x_j,x_k)\leq r_j-r+r_k$$
for all $j,k\in J$, there exists $z\in X$ such that
$$\al(x_j,z)\leq r_j$$
for all $j\in J$. Let us elaborate this description as follows:
\begin{itemize}
\item Each pair $(x_j,r_j)$ may be viewed as a closed ball $B(x_j,r_j)$ of center $x_j$ and radius $r_j$. Since, as a partial metric space, the distance $\al(x_j,x_j)$ of $x_j$ to itself may not be zero, the inequality $\al(x_j,x_j)\leq r_j$ says that the radius of the ball $B(x_j,r_j)$ cannot be less than the self-distance of its center.
\item The requirement $r_j\in[r,\infty]$ means that there is a closed ball $B(x_j,r)$ of fixed radius $r$ inside every ball $B(x_j,r_j)$. Moreover, the inequality $\al(x_j,x_k)\leq r_j-r+r_k$ says that any two balls $B(x_j,r_j)$ and $B(x_k,r_k)$ ``strongly intersect'', which intuitively means that the straight line connecting their centers has at least a segment of length $r$ lying inside the intersection of $B(x_j,r_j)$ and $B(x_k,r_k)$.
\item The existence of $z\in X$ such that $\al(x_j,z)\leq r_j$ for all $j\in J$ means that the intersection $\bigcap\limits_{j\in J}B(x_j,r)$ of all balls is non-empty. 
\end{itemize}
\end{rem}

As a special case of Theorem \ref{hyperconvex=injective}, hyperconvex partial metric spaces are precisely injective partial metric spaces:

\begin{cor} \label{parmet-hyperconvex=injective}
Let $(X,\al)$ be a partial metric space. Then the following statements are equivalent:
\begin{enumerate}[label={\rm(\roman*)}]
\item \label{parmet-hyperconvex=injective:hyperconvex} $(X,\al)$ is hyperconvex.
\item \label{parmet-hyperconvex=injective:AR} $(X,\al)$ is an absolute retract in $\ParMet$.
\item \label{parmet-hyperconvex=injective:injective} $(X,\al)$ is an injective object in $\ParMet$.
\end{enumerate}
\end{cor}

Now let us define the tight span $(\sT(X,\al),\si)$ of a partial metric space $(X,\al)$. As a preparation, we point out that implications in the quantaloid $\sD[0,\infty]$ \cite{Tao2014} are given by
$$w\lda u=q\vee r\vee(w+q-u)\quad\text{and}\quad v\rda w=p\vee q\vee(w+q-v)$$
for all $u:p\to q$, $v:q\to r$, $w:p\to r$ in $\sD[0,\infty]$. Hence, it follows from Remark \ref{TX-def-QRel} that $\sT(X,\al)$ consists of maps $\mu:X\to[r,\infty]$ $(r\in[0,\infty])$ satisfying
$$\al(x,x)\leq\mu(x)=r\vee\al(x,x)\vee\sup\limits_{y\in Y}(\al(x,y)+r-\mu(y))$$
for all $x\in X$. The partial metric $\si$ on $\sT(X,\al)$ is given by
$$\si(\mu,\lam)=r\vee s\vee\sup\limits_{x\in X}(\lam(x)+r-\mu(x))$$
for all $\mu:X\to[r,\infty]$, $\lam:X\to[s,\infty]$ in $\sT(X,\al)$. As an immediate consequence of Theorem \ref{injective-hull}:

\begin{cor} \label{tight-span} 
The tight span $(\sT(X,\al),\si)$ is the injective hull of a partial metric space $(X,\al)$ in $\ParMet$.
\end{cor}

Finally, for the essentiality of \emph{isometric maps} (i.e., fully faithful $\sD[0,\infty]$-functors) between partial metric spaces, we translate Theorem \ref{f-essential} to the following:

\begin{cor} \label{f-essential-metric}
Let $f:(X,\al)\to(Y,\be)$ be an isometric map between partial metric spaces. Then $f$ is essential in $\ParMet$ if, and only if, $f$ is dense in the sense that
$$\be(y,y')=\be(y,y)\vee\be(y',y')\vee\sup\limits_{x\in X}(\be(fx,y')+\be(y,y)-\be(fx,y))$$
for all $y,y'\in Y$.
\end{cor}

\section*{Acknowledgements}

The authors would like thank constructive comments received from the anonymous referee, as well as helpful discussions with Professor Hongliang Lai and Professor Dexue Zhang.



\end{document}